\documentclass[12pt]{article}

\usepackage[colorlinks,allcolors=blue]{hyperref}

\usepackage[margin=1in]{geometry}
\usepackage{amsfonts, amsmath, amssymb, amsthm} 
\usepackage[none]{hyphenat}
\usepackage{fancyhdr} 
\usepackage[utf8]{inputenc}
\usepackage{mathtools}

\makeatletter
\def\legendre@dash#1#2{\hb@xt@#1{%
  \kern-#2\p@
  \cleaders\hbox{\kern.5\p@
    \vrule\@height.2\p@\@depth.2\p@\@width\p@
    \kern.5\p@}\hfil
  \kern-#2\p@
  }}
\def\@legendre#1#2#3#4#5{\mathopen{}\left(
  \sbox\z@{$\genfrac{}{}{0pt}{#1}{#3#4}{#3#5}$}%
  \dimen@=\wd\z@
  \kern-\p@\vcenter{\box0}\kern-\dimen@\vcenter{\legendre@dash\dimen@{#2}}\kern-\p@
  \right)\mathclose{}}
\newcommand\legendre[2]{\mathchoice
  {\@legendre{0}{1}{}{#1}{#2}}
  {\@legendre{1}{.5}{\vphantom{1}}{#1}{#2}}
  {\@legendre{2}{0}{\vphantom{1}}{#1}{#2}}
  {\@legendre{3}{0}{\vphantom{1}}{#1}{#2}}
}
\def\dlegendre{\@legendre{0}{1}{}}
\def\tlegendre{\@legendre{1}{0.5}{\vphantom{1}}}
\makeatother

\usepackage{booktabs, siunitx,caption}

\newtheorem{theorem}{Theorem}

\title{Relaxing the size constraints on Proth's criterion}
\author{Tejas R. Rao}
\date{}

\begin{document}
\maketitle

\begin{center}
\small \textbf{Abstract.} We add one condition to Proth's theorem to extend its applicability to $N=k2^n+1$ where $2^n>N^{1/3}$ as opposed to the former constraint of $2^n>k$. This additional condition adds barely any complexity or time to the test and can furthermore be calculated concurrently. Furthermore, it maintains the biconditionality of Proth's theorem and thus makes it readily applicable. A note on an extension of Brillhart, Lehmer, and Selfridge's primality test is also made. 
\end{center}

\setlength{\parindent}{0cm}

The famous Proth primality test is an adaption of Pocklington's criterion and has been the subject of dedicated computation (\hyperlink{1}{1}). It states that $N=k2^n+1$, where $k$ is odd and $2^n>k$ is prime if and only if 
\begin{center}
$a^{\frac{N-1}{2}}\equiv 1\mod N$,
\end{center}
for all $a$ where $\legendre{a}{N}=-1$. However, the size constraints require that $2^n>k$, or about that $2^n>\sqrt{N}$. Another famous extension of Pocklington's criterion is as follows. From Brillhart, Lehmer, and Selfridge, we know that, for $N=mp+1$, where $p$ is prime and $m$ is a positive integer, if 
\begin{align*}
2p+1&>\sqrt{N},\\
a^{N-1}&\equiv 1\mod N,\\
a^{m}&\not\equiv 1\mod N,
\end{align*}
then $N$ is prime (\hyperlink{5}{5}). Note that this criterion is significantly weaker because it is not biconditional. Proth's test is deterministic for any chosen $N$ because $\legendre{a}{N}$ may be easily calculated using quadratic reciprocity (\hyperlink{2}{2}, \hyperlink{4}{4}). However, the test of Brillhart, Lehmer, and Selfridge requires random chance in choosing an $a$, and additionally does not prove $N$ is composite. In this paper, we relax the size constraints on $2^n$ in Proth's test by utilizing the techniques found in (\hyperlink{1}{1}). 
\vspace{5mm}\\
We begin with the following proposition, utilizing a similar method of proof as Brillhart, Lehmer, and Selfridge do (\hyperlink{5}{5}). 
\begin{theorem}
For $N=mp^z+1$, if 
\begin{align*}
2p^z+1&>\sqrt{N},\\
a^{N-1}&\equiv 1\mod N,\\
a^{mp^{z-1}}&\not\equiv 1\mod N,
\end{align*}
then $N$ is prime. 
\end{theorem}
If $a^{mp^z}\equiv 1\mod N$ but $a^{mp^{z-1}}\not\equiv 1\mod N$, then $p^z|ord_N(a)$. If $p^z|ord_N(a)$, and since $ord_N(a)=lcm[ord_{n_i^{z_i}}(a)]_i$, where $n_i^{z_i}$ represent the highest powers of each prime factor of $N$, at least one prime factor $n_i$ must satisfy $p^z|ord_{n_i^{z_i}}(a)$. But since, from (\hyperlink{3}{3}), $ord_{n_i^{z_i}}(a)=n_i^{q}ord_{n_i}(a)$ for some $q$, and since $q\nmid N$, we know $p^z|ord_{n_i}(a)$. Therefore, since $n_i$ is prime, we can write 
\begin{center}
$p^z|ord_{n_i}(a)|\phi(n_i)=n_i-1$. 
\end{center}
Alternatively, $n_i\equiv 1\mod p^z$. So, $N/n_i\equiv 1\mod p^z$. But since $n_i\geq 2p^z+1>\sqrt{N}$, we know $N/n_i<\sqrt{N}<2p^z+1$, and thus $N/n_i$ is precisely $1$ and $N=n_i$ is prime. 
\vspace{5mm}\\
This extension is valid, but not as strong as that of Proth because, again, it is not biconditional. However, utilizing quadratic reciprocity and techniques from (\hyperlink{1}{1}), we can relax the size criterion on Proth's primality test. Recall that for all primes, 
\begin{center}
$a^{\frac{N-1}{2}}\equiv\legendre{a}{N}=-1\mod N$
\end{center}
for all $a$ where $\legendre{a}{N}=-1$. 
Combine this with the following theorem. 
\begin{theorem}
For $N=k2^n+1$, where $k$ is odd and $2^n>\sqrt[3]{N}$, $N$ is either prime or semiprime if
\begin{center}
$a^{\frac{N-1}{2}}\equiv -1\mod N$.
\end{center}
\end{theorem}
By the conditions, $a^{N-1}\equiv 1\mod N$, and thus $2^n|ord_N(a)$. 
Additionally, since it is congruent to $-1$, we know that, for each prime factor $n_i$ of $N$, $a^{\frac{N-1}{2}}-1\equiv -1-1=-2\mod n_i$. Since $N$ is odd, all prime factors are greater than $2$ and thus
\begin{center}
$gcd(a^{\frac{N-1}{2}}-1, N)=1$.
\end{center} 
Because of this we know that $2^n|ord_{n_i}(a)$ for all prime factors $n_i$ of $N$. Thus, 
\begin{center}
$2^n|ord_{n_i}(a)|\phi(n_i)=n_i-1$. 
\end{center}
Therefore, $n_i\equiv 1\mod 2^n$ for all prime factors $n_i$. But also $n_i\geq 2^n+1>\sqrt[3]{N}$. So $N/n_1n_2<N^{1/3}<2^n+1$. But since $N/n_1n_2\equiv 1\mod 2^n$, so $N/n_1n_2$ is precisely $1$ and $N$ is thus either prime or semiprime. 
\vspace{5mm}\\
Combining the two aforementioned statements, we know that for a prime $N=k2^n+1$ with $2^n+1>\sqrt[3]{N}$, $a^{\frac{N-1}{2}}\equiv -1\mod N$. In the other direction, if $a^{\frac{N-1}{2}}\equiv -1\mod N$, then $N$ is either prime or semiprime. We can easily determine whether $N$ is semiprime given the constraints. Specifically, we utilize the fact that if $q<w$, then $q\mod w = q$. If $N$ is semiprime and satisfies the aforementioned conditions, then we may write $N=(2^nu+1)(2^nv+1)=2^{n}(2^nuv+u+v)+1$. Since $2^n2^nuv<N$, we know $u+v<uv<N^{1/3}<2^n$ (if $u+v\geq uv$, then because of the small values of $u,v$, the inequality less than $2^n$ will still hold). Therefore, we can find $u+v$ by calculating $k\mod 2^n$ and $uv$ by calculating $\frac{k-k\mod 2^n}{2^n}$. Solving this system of equations, we can find $u,v$ and thus the prime factors $2^nu+1$ and $2^nv+1$. Thus, if and only if the solution is not an integer, $N$ is prime. This means we have the following conditions after solving the system. 
\begin{theorem}
For $N=k2^n+1$, $k$ odd, $\legendre{a}{N}=-1$, and $2^n>N^{1/3}$, $N$ is prime if and only if 
\begin{align*}
a^{\frac{N-1}{2}}&\equiv -1\mod N,\\
u&\not\in\mathbb{N},
\end{align*}
where 
\begin{center}
$u=\frac{1}{2}(k\mod 2^n-\sqrt{(k\mod 2^n)^2-4(\frac{k-k\mod 2^n}{2^n})}).$
\end{center}
\end{theorem}
Furthermore, if $u\in\mathbb{N}$, then $N$ is semiprime and its factors are given by $2^nu+1$ and $2^nv+1$ as defined above. The second condition only requires a few modular multiplications to calculate and thus does not add significant complexity to the test. Additionally, regardless of the outcome of the second condition, at least one and up to two primes will be discerned ($N$ or its factors). If $2^n>N^{1/2}$, the second condition may be removed.

\newcommand{\noopsort}[1]{} \newcommand{\printfirst}[2]{#1}
  \newcommand{\singleletter}[1]{#1} \newcommand{\switchargs}[2]{#2#1}

\end{document}